\documentclass[11pt,reqno]{amsart}
\usepackage{layout}
\usepackage{latexsym,amsmath,amssymb,amsfonts,amsthm,verbatim}
\usepackage[utf8]{inputenc}
\usepackage{color}

\let \ssection=\section
\renewcommand{\section}{\setcounter{equation}{0}\ssection}
\def\lv{\left|}
\def\rv{\right|}
\def\rvd{\right|^2}
\def\undeux{{1\over2}}
\def\bknlm{\beta^{k,N}_\lm}
\def\gknlm{\gamma^{k,N}_\lm}
\def\txnlm{\lp t,\xnlm(t)\rp}
\newtheorem{Th}{Theorem}[section]
\newtheorem{theorem}[Th]{Theorem}

\newtheorem{lemma}[Th]{Lemma}
	\theoremstyle{definition}

	\theoremstyle{remark}

\def\lc{\left[}
\def\rc{\right]}
\def\H{{L^2(0,1)}}
\def\undeux{{{1\over2}}}
\def\A{{{1\over2}{d^2\over dx^2}}}

\begin{document}
\title{Weak error in negative Sobolev spaces for the stochastic heat equation}
\author{Omar Aboura}
 \email{omar.aboura@malix.univ-paris1.fr}
\address{ 
  SAMM, EA 4543,
 Universit\'e Paris 1 Panth\'eon Sorbonne,
90 Rue de Tolbiac, 75634 Paris Cedex France }
\begin{abstract}
In this paper, we make another step in the study of weak error of the
stochastic heat equation by considering norms as functional.
\end{abstract}
\maketitle

\section{Introduction}
Let $(\Omega, {\mathcal F}, P)$ a probability space and $T>0$ a fixed time.
$(W(t))_{t\geq0}$ will be a cylindrical Brownian motion on $\H$.
We consider the stochastic heat equation, written in abstract form in
$\H$: $X(0)=0$, for all $t\in[0,T]$ $X(t,0)=X(t,1)=0$ and
\begin{equation}dX(t)=\A X(t)dt+dW(t).\label{eq-X}
\end{equation}
It is well know that this equation admits 
a unique weak solution (from the analytical point of view).

Let $N\in\mathbb N^*$ and $h:=T/N$.
Consider $(t_k)_{0\leq k\leq N}$ the uniform subdivision of $[0,T]$
defined by $t_k:=kh$.
We consider the implicit Euler scheme defined as follow:
\def\XN{{X^N}}
\def\ltkk{{(t_{k+1})}}
\def\ltk{{(t_k)}}
\def\DWkk{\Delta W(k+1)}
\begin{equation}\XN\ltkk=\XN\ltk+h\A\XN\ltkk+\DWkk,\label{eq-XN}
\end{equation}
where $\DWkk=W\ltkk-W\ltk$.

\def\R{{\mathbb R}}
\def\lp{\left(}
\def\rp{\right)}
Let $f:\H\rightarrow\R$ be a functional.
The stong error is the study of $E\lv\XN(T)-X(T)\rv^2_{L^2(0,1)}$.
The weak error is the study of $\lv Ef\lp\XN(T)\rp-Ef\lp X(T)\rp\rv$
with respect to the time mesh $h$.

In \cite{Debussche}, A.~Debussche
considers a more general stochastic equation and a more general functional
than the one considered here. He obtains a weak error of order $1/2$,
which is the double of that proved by \cite{PJ} for the strong speed of
convergence.
The novelty of this paper his to prove that for the square of the norm the weak
error his better than $1/2$ in negative Sobolev spaces.
\section{Preliminaries and main result}
\subsection*{Notations}
We collect here some of the notations used through the paper.
 $<.,.>_{\H}$ is the inner product in $L^2(0,1)$,
 $H^1_0(0,1)$ is the Sobolev space of functions $f$ in $\H$
vanishing in 0 and 1 with first derivatives in $\H$,
$H^2(0,1)$ is the Sobolev space of functions $f$ in $\H$ with first and second derivatives in $\H$.
Finally, for $m=1,2,\dots$, let $(e_m(x)=\sqrt{2}\sin(m\pi x)$ and
$\lambda_m=\undeux(\pi m)^2$ denote the eigenfunction and eigenvalues of $-\Delta$
with Dirichlet boundary conditions on $(0,1)$.

An $\H$-valued stochastic process $\lp X(t)\rp_{t\in[0,T]}$
is said to be a solution of \eqref{eq-X} if:
$X(0)=0$ and for all $g\in H^1_0(0,1)\cap H^2(0,1)$ we have
$$<X(t),g>_\H=\int_0^t<X(s),\A g>_\H ds+<W(t),g>_\H.$$
It is well know that \eqref{eq-X} admits a unique solution: see \cite{DaPrato}.
Then $(e_m)_{m\geq1}$ is a complete orthonormal basis of $\H$.
\def\lm{{\lambda_m}}
\def\xlm{{X_\lm}}
If we denote by $\lm:=\undeux(\pi m)^2$, $W_\lm(t):=\left<W(t),e_m\right>_H$
and $\xlm(t)$ denote the solution of the evolution equation:
$X_\lm(0)=0$ and for $t>0$:
$$d\xlm(t) = -\lm \xlm(t)dt+dW_\lm(t).$$
Then the processes $\lp\xlm(.)\rp_{m\geq1}$ are independent and
$X(t)=\sum_{m\geq1}\xlm(t)e_m$ for all $t\geq0$.

A sequence of $\H$-valued $\lp\XN\ltk\rp_{k=0,\dots,N}$ is said to be a
solution of \eqref{eq-XN} if:
$\XN(t_0)=0$ and for all $k=0,\dots,N-1$ and for all
$g\in H^1_0(0,1)\cap H^2(0,1)$ we have
\def\xn{{X^N}}
\begin{align*}
<\xn\ltkk,g>_\H =& <\xn\ltk,g>_\H +h<\xn\ltkk,\A g>_\H\\
	&+<\DWkk,g>_\H.
\end{align*}
It is well know that \eqref{eq-XN} has a unique solution and there exists
a constant $C>0$, independent of $N$, such that
$E\lv\xn(T)-X(T)\rvd_\H\leq Ch^{\undeux}$ where $h=T/N$.
\def\xnlm{{X^N_\lm}}
Now if we denote by $\lp\xnlm\ltk\rp_{k=0,\dots,N}$ the solution of:
$\xnlm(t_0)=0$ and for $k=0,\dots,N-1$
$$\xnlm\ltkk = \xnlm\ltk -\lm h\xnlm\ltkk +W_\lm(k+1).$$
The random vectors $(\xnlm(t_k), k=0,\dots,N)_{m=1,2,\dots}$ are independent and
$\xn\ltk=\sum_{m\geq1}\xnlm\ltk e_m$.

\def\Hp{{H^{-p}}}
\def\summ{\sum_{m\geq1}}
Let $p\geq0$; we define the spaces $\Hp$ as the completion of $\H$
for the topology induced by the norm 
\def\lmp{\lambda^{-p}_m}
$\lv u\rvd_\Hp:=\summ\lmp<u,e_m>_H^2$.
The following theorem improves the speed of convergence of $X^N$ to $X$
for negative Sobolev spaces.
\begin{theorem} \label{th1}
Suppose that $h<1$ and let $p\in[0,\undeux)$.
There exists a constant $C>0$, independent of $N$, such that
$$\lv E\lv\xn(T)\rvd_\Hp-E\lv X(T)\rvd_\Hp\rv\leq C h^{p+\undeux}.$$
\end{theorem}

\section{Proof of the theorem \ref{th1}}
The proof of the theorem will be done in several steps.
First we recall the weak error of the Ornstein-Uhlenbeck process.
Secondly we prove some technical lemmas.
Then we decompose the weak error and analyse each term of
these decomposition.

\subsection{Weak error of the Ornstein-Uhlenbeck process}
\def\wl{W_{\lambda}}
\def\xl{X_{\lambda}}
Let $\lambda>0$, $(\wl(t))_{t\geq0}$ be a one dimensional Brownian motion
and $(\xl(t))_{t\geq0}$ be the Ornstein Uhlenbeck process solution of the following stochastic 
differential equation:
$\xl(0)=x\in\R$ and
\begin{equation}d\xl(t)=-\lambda\xl(t)dt+d\wl(t).\label{eq-ou}
\end{equation}

In this step, we study two properties associated with this process:
the Kolmogorov equation and the implicit Euler scheme.

Let $\lp X^{t,x}_{\lambda}(s)\rp_{t\leq s\leq T}$ be the solution of \eqref{eq-ou}
starting from $x$ at time $t$.
It is well know that $X^{t,x}_{\lambda}(T)$ is a normal random variable:
$$X^{t,x}_{\lambda}(T) \sim
\mathcal{N}\lp e^{-\lambda(T-t)}x,{1-e^{-2\lambda(T-t)}\over 2\lambda}\rp.$$
For $t\in[0,T]$ and $x\in\R$ set $u_{\lambda}(t,x):=E\lv X^{t,x}_{\lambda}(T)\rvd$.
Then $u_{\lambda}$ is the solution of the following partial differential equation,
called Kolmogorov equation:
for all $x\in\R$, $u_{\lambda}(T,x)=\lv x\rvd$
and for all $(t,x)\in[0,T)\times\R$
\def\dt{{\partial\over\partial t}}
\def\dx{{\partial\over\partial x}}
\def\dxx{{\partial^2\over\partial x^2}}
\begin{equation}
-\dt u(t,x)=\undeux\dxx u(t,x)-\lambda x\dx u(t,x).\label{eq-kol}
\end{equation}
Since $X^{t,x}_{\lambda}(T)$ has a normal law, we can write $u_{\lambda}$
explicitely:
\begin{equation}
u_{\lambda}(t,x)={1-e^{-2\lambda(T-t)}\over2\lambda}
 +e^{-2\lambda(T-t)}x^2.\label{eq-u}
\end{equation}
With this expression we see that $u_{\lambda}\in C^{1,2}([0,T]\times\R)$
and we have the following derivatives:
\def\ultx{u_{\lambda}(t,x)}
\def\eltt{e^{-2\lambda(T-t)}}
\def\dtx{{\partial^2\over\partial t\partial x}}
\begin{align}
\dx\ultx=&2\eltt x,\label{eq-ux}\\
\dxx\ultx=&2\eltt,\label{eq-uxx}\\
\dt\ultx=&-\eltt+2\lambda\eltt x^2,\label{eq-ut}\\
\dtx\ultx = &4\lambda\eltt x.\label{eq-utx}
\end{align}

\def\xnl{X^N_{\lambda}}
\def\Wl{{W_\lambda}}
\def\ulh{{1\over1+\lambda h}}
\def\dwl{\Delta W_{\lambda}}
The implicit Euler scheme for the Ornstein-Uhlenbeck equation \eqref{eq-ou}
starting from 0 at time $t_0$, is defined as follow:
$\xnl(t_0)=0$ and for $k=0,\dots,N-1$
\begin{equation}\xnl\ltkk=\xnl\ltk-\lambda h\xnl\ltkk
+\Delta W_{\lambda}(k+1),\label{eq-xnl}
\end{equation}
where $\Delta W_{\lambda}(k+1)=\Wl\ltkk-\Wl\ltk$.
Since we have the following equation
\begin{equation}
\xnl\ltkk=\ulh\xnl\ltk+\ulh\dwl(k+1),\label{eq-xnld}
\end{equation}
we see that the scheme is well defined.
\begin{lemma} \label{lem-xsum}
For $k=1,\dots,N$ we have
$\xnl\ltk=\sum_{j=0}^{k-1}{\dwl(k-j)\over(1+\lambda h)^{j+1}}.$
\end{lemma}
\proof We proceed by induction.
If $k=1$, we have $\xnl(t_1)=\ulh\dwl(1)$.
Suppose the result true until $k$.
Using \eqref{eq-xnld}, we have
\begin{align*}
\xnl\ltkk=&\sum_{j=0}^{k-1}{\dwl(k-j)\over(1+\lambda h)^{j+2}}
	+\ulh\dwl(k+1)\\
=&\sum_{l=1}^k{\dwl(k+1-l)\over(1+\lambda h)^{l+1}}
	+{1\over(1+\lambda h)^{0+1}}\dwl(k+1-0),
\end{align*}
which concludes the proof.
\endproof
\begin{lemma}\label{lem-EX}
For all $k=0,\dots,N$, we have the following bound
$E\lv\xnl\ltk\rvd \leq {1\over2\lambda}.$
\end{lemma}
\proof
Using the independence of the increments of the Brownian motion and  
Lemma \ref{lem-xsum}, we have
\begin{align*}
E\lv\xnl\ltk\rvd =& \sum_{j=0}^{k-1} {1\over(1+\lambda h)^{2(j+1)}}
			E\lv\dwl(k-j)\rvd
		=  h\sum_{j=0}^{k-1} {1\over(1+\lambda h)^{2(j+1)}}.
\end{align*}
\def\lh{\lambda h}
Let $a:=1/(1+\lh)^2$; we deduce that $E\lv\xnl\ltk\rvd = ha {1-a^k\over1-a}$.
Simple computations yield $ha/(1-a) = 1/(2\lambda+\lambda^2h)$,
which implies
$$E\lv\xnl\ltk\rvd = {1\over2\lambda+\lambda^2h}
			\lp 1-{1\over(1+\lh)^{2k}}\rp.$$
This concludes the proof.
\endproof
For $t\geq 0$, we denote $\mathcal{F}^{\lambda}_t :=
\sigma\lp W_\lambda(s), s\leq t\rp$ and $D^{1,2}_\lambda$
the Malliavin Sobolev space with respect to $W_\lambda$.
\begin{lemma} \label{lem-mal}
For all $k=1,\dots,N$, we have
$\xnl\ltk\in D_\lambda^{1,2}\cap L^2\lp\mathcal{F}^\lambda_{t_k}\rp.$
\end{lemma}
\proof This is a consequence of Lemma \ref{lem-xsum}, the fact that 
$L^2\lp\mathcal{F}^\lambda_{t_k}\rp$ and $D_\lambda^{1,2}$ are linear space and for all
$j=0,\dots,k-1$, $\dwl(k-j)\in D_\lambda^{1,2}\cap L^2\lp\mathcal{F}^\lambda_{t_k}\rp$.
\endproof

As usual in the study of weak error, we need to use a continuous process that
interpolates the Euler scheme.
The interpolation process that we use was introduced in \cite{Aboura}.
We recall its construction and prove some of its properties.

Let $k\in\{0,\dots,N-1\}$ be fixed.
In order to interpolate the scheme between the points $\lp t_k,\xnl\ltk\rp$
and $\lp t_{k+1},\xnl\ltkk\rp$, we define the process as follows:
for $t\in[t_k,t_{k+1}]$, set
\begin{equation}
\xnl(t) := \xnl\ltk-\lambda E\lp\xnl\ltkk|\mathcal{F}_t\rp(t-t_k)
	+W_{\lambda}(t)-W_{\lambda}\ltk. \label{eq-xns}
\end{equation}
In the sequel, we will use the following processes:
for $t\in[t_k,t_{k+1}]$
\def\F{\mathcal{F}}
\def\tkk{{t_{k+1}}}
\def\bknl{\beta^{k,N}_{\lambda}}
\def\zknl{z^{k,N}_{\lambda}}
\def\gknl{\gamma^{k,N}_{\lambda}}
\begin{align}
\bknl(t) :=& -\lambda E\lp\xnl\ltkk|\F_t\rp, \label{eq-beta}\\
\zknl(t) :=& -\lambda E\lp D_t\xnl\ltkk|\F_t\rp, \label{eq-z}\\
\gknl(t) :=& 1+(t-t_k)\zknl(t). \label{eq-gamma}
\end{align}
The next lemma relates the above processes.
\def\wl{W_{\lambda}}
\begin{lemma}\label{lem-bzX}
Let $k=0,\dots,N-1$. 
For $t\in[0,T]$, we have
\begin{align*}
d\bknl(t) =& \zknl(t)d\wl(t),\quad
\zknl(t) = -{\lambda\over1+\lambda h},\\
\gknl(t) =& 1 - (t-t_k){\lambda\over1+\lambda h},\quad
d\xnl(t) = \bknl(t)dt + \gknl(t)d\wl(t).
\end{align*}
\end{lemma}
\proof Using the Clark-Ocone formula and Lemma \ref{lem-mal}, we have
$$\xnl\ltkk = E\lp\xnl\ltkk|\F_t\rp 
+ \int_t^\tkk E\lp D_s\xnl\ltkk|\F_s\rp d\wl(s).$$
Multiplying by $(-\lambda)$, we deduce
$$-\lambda\xnl\ltkk = \bknl(t)+\int_t^\tkk\zknl(s) d\wl(s),$$
which gives the first identity.
Applying the Malliavin derivative to \eqref{eq-xnld},
we have for $s\in[t_k,\tkk]$
$D_s\xnl\ltkk = \ulh$.
Multiplying by $(-\lambda)$, we deduce the second and third equalities.

Finaly, It\^o's formula gives us
$$d\lp(t-t_k)\bknl(t)\rp= (t-t_k)\zknl(t)d\wl(t)+\bknl(t)dt,$$
which concludes the proof.
\endproof
\begin{lemma}\label{lem-bx}
Let $k\in\{0,\dots,N-1\}$.
For any $s\in[t_k,\tkk]$, we have
\begin{align*}
E\lv\bknl(s)\rvd \leq& 2\lambda,\quad
E\lv\xnl(s)\rvd \leq {1\over2\lambda} + h,\quad
E\bknl(s)\xnl(s) \leq 1.
\end{align*}
\end{lemma}
\proof 
Applying the conditionnal expectation with respect to $\F_s$
on both sides of \eqref{eq-xnld} for $s\in[t_k,t_{k+1})$ we have
$$E\lp\xnl\ltkk|\F_s\rp = \ulh\lc\xnl\ltk + (\wl(s)-\wl\ltk)\rc.$$
Multiplying by $(-\lambda)$ and using \eqref{eq-beta}, we obtain
\begin{equation}
\bknl(s) = -{\lambda\over 1+\lambda h}\xnl\ltk
	-{\lambda\over1+\lambda h}\lp\wl(s)-\wl\ltk\rp. \label{eq-betad}
\end{equation}
The independence of $\F_{t_k}$ and $\wl(s)-\wl\ltk$ yields
$$E\lv\bknl(s)\rvd = {\lambda^2\over(1+\lambda h)^2} E\lv\xnl\ltk\rvd
	+{\lambda^2\over(1+\lambda h)^2}(s-t_k).$$
Using Lemma \ref{lem-EX}, we deduce
$$E\lv\bknl(s)\rvd \leq {\lambda\over2(1+\lambda h)^2}
	+ {\lambda^2h\over(1+\lambda h)^2},$$
which proves the first upper estimate.

Using \eqref{eq-xns} and \eqref{eq-betad}, we have for $s\in[t_k,t_{k+1}]$
\def\AA{\lp1-{\lambda(s-t_k)\over1+\lambda h}\rp}
\begin{equation}
\xnl(s) = \AA\lc\xnl\ltk + (\wl(s)-\wl\ltk)\rc. \label{eq-xnsd}
\end{equation}
Taking the expectation of the square and using the independence
of $\F_{t_k}$ and $\wl(s)-\wl\ltk$, we have
\begin{align*}
E\lv\xnl(s)\rvd =& \AA^2\lc E\lv\xnl\ltk\rvd + (s-t_k) \rc
	\leq  E\lv\xnl\ltk\rvd + h
	\leq {1\over 2\lambda} + h,
\end{align*}
where the last upper estimates follows from Lemma \ref{lem-EX}.

Multiplying \eqref{eq-betad} and \eqref{eq-xnsd}, taking expectation we obtain
\def\BB{{-\lambda\over1+\lambda h}}
\def\sqa{\sqrt{{q\over\alpha}}}
\def\csqa{\left[\sqa\right]}
\def\mp{{-p}}
\def\xnlm{X^N_\lm}
\def\xlm{X_\lm}
\def\dnkm{\delta^N(k,m)}
\def\lb{\left\lbrace}
\def\rb{\right\rbrace}
\def\tk{{t_k}}
$$E\lp\xnl(s)\bknl(s)\rp = \BB\AA\lc E\lv\xnl\ltk\rvd + (s-t_k)\rc.$$
Using Lemma \ref{lem-EX}, we deduce
$$\lv E\lp\xnl(s)\bknl(s)\rp\rv \leq {\lambda\over1+\lambda h}{1\over2\lambda}
	+ {\lambda h\over1+\lambda h}.$$
This concludes the proof.
\endproof

\subsection{Some useful analytical lemmas}
We at first give a precise upper bound of a series defined in terms of the
eigenvalues of the Laplace operator with Dirichlet boundary conditions.

\begin{lemma}\label{lem-a1}
Let $p\in[0,\undeux)$. There exists a constant $C>0$, 
such that for all $\alpha>0$, we have
$$\summ\lmp e^{-2\lm\alpha}\leq C\alpha^{p-\undeux}$$
\end{lemma}
\proof
The function $(x\in\R_+\mapsto x^{-2p}e^{-2x^2\alpha})$ is decreasing.
So by comparaison, we obtain
\begin{align*}
\summ m^{-2p}e^{-2m^2\alpha} \leq \int_0^{\infty}x^{-2p} e^{-2x^2\alpha}dx
	\leq \alpha^{p-\undeux}\int_0^\infty y^{-2p}e^{-2y^2}dy
	=  C \alpha^{p-\undeux}.
\end{align*}
Since $\lm = \undeux(\pi m)^2$, we deduce the desired upper estimate.
\endproof
\begin{lemma}\label{lem-ad}
Let $q>0$.
There exists a constant $C>0$, such that for all $\alpha>0$
$$\summ\lambda^q_me^{-\lm\alpha} \leq C\lp1+{1\over\alpha^{q+\undeux}}\rp.$$
\end{lemma}
\proof Let $f(x) = x^{2q}e^{-x^2\alpha}$.
His derivatives is given by $f'(x) = 2x^{2q-1}e^{-x^2\alpha}(q-\alpha x^2)$.
\textit{Case 1:} $\alpha>q/4$.
Then $f$ is decreasing on $[2,\infty)$ and a standard comparaison argument yields
\begin{align*}
\summ m^{2q}e^{-m^2\alpha} \leq & e^{-\alpha} + 4^qe^{-4\alpha}
	+\sum_{m\geq3}\int_{m-1}^{m}x^{2q}e^{-x^2\alpha}dx \\
	\leq & C + \int_0^\infty x^{2q}e^{-x^2\alpha} dx \\
	\leq & C + \alpha^{-q-\undeux}\int_0^\infty y^{2q}e^{-y^2}dy \\
	\leq & C (1+\alpha^{-q-\undeux}).
\end{align*}
\textit{Case 2:} $\alpha \leq q/4$.
The function $f$ is increasing on $[0,\sqrt{q/\alpha}]$.
So for each $m=1,\dots,[\sqa]-1$, we have
$$m^{2q}e^{-m^2\alpha} \leq \int_m^{m+1} x^{2q} e^{-x^2\alpha}dx.$$
On the interval $[\sqa,\infty)$, $f$ is decreasing.
So for each integer $m\geq[\sqa]+2$, we have
$$m^{2q}e^{-m^2\alpha} \leq \int_{m-1}^mx^{2q}e^{-x^2\alpha}dx.$$
The above upper estimates yield
\begin{align*}
\summ m^{2q}e^{-m^2\alpha} \leq &
	\sum_{m\leq[\sqa]-1}\int_m^{m+1} x^{2q} e^{-x^2\alpha}dx
	+ \sum_{m\geq[\sqa]+2}\int_{m-1}^m x^{2q} e^{-x^2\alpha}dx\\
	& + \sum_{m\in\{[\sqa],[\sqa]+1\}}m^{2q}e^{-m^2\alpha}\\
	\leq &\int_0^\infty x^{2q} e^{-x^2\alpha}dx
	+ \sum_{m\in\{[\sqa],[\sqa]+1\}}m^{2q}e^{-m^2\alpha}\\
	\leq &C\alpha^{-q-\undeux}
	+ \sum_{m\in\{[\sqa],[\sqa]+1\}}m^{2q}e^{-m^2\alpha}
\end{align*}
Now we study each term of the sum in the right hand side.
Since $q\geq\alpha$, we have
\begin{align*}
\csqa^{2q}e^{-\csqa^2\alpha} \leq & \lp {q\over \alpha}\rp^q
	\leq  \lp{q\over\alpha}\rp^{q+\undeux}
	\leq  C \alpha^{-q-\undeux}.
\end{align*}
For the second term, we remark that since $q\geq\alpha$ $\csqa+1\leq2\csqa\leq2\sqa$.
This implies
\begin{align*}
\lp\csqa+1\rp^{2q}e^{-\lp\csqa+1\rp^{2}\alpha} \leq \lp 2\sqa\rp^{2q} 
	\leq C \alpha^{-q-\undeux}.
\end{align*}
Therefore, in both cases we obtain
$$\summ m^{2q}e^{-m^2\alpha} \leq C\lp 1+{1\over\alpha^{q+\undeux}}\rp.$$
Since $\lm=\undeux(\pi m)^2$, the proof is complete.
\endproof
\begin{lemma}\label{lem-at}
Let $p\in[0,\undeux)$ and $n\in\mathbb{N}^*$.
Let $\lp v(k,m)\rp_{(k,m)\in\{0,\dots,N-2\}\times\mathbb{N}^*}$ be a sequence
such that for all $k\in\{0,\dots,N-2\}$ and $m\geq 1$, we have
$$0\leq v(k,m) \leq \lambda^{n-p}_mh^{n+1}e^{-2\lm(T-\tkk)}.$$
Then, there exists a constant $C>0$, independent of $N$, such that
$$\summ\sum_{k=0}^{N-2} v(k,m)\leq C h^{p+\undeux}.$$
\end{lemma}
\proof
First we remark that $T-\tkk = h(N-k-1)$.
Using Lemma \ref{lem-ad}, we deduce the existence of $C$ 
depending on $n$ and $p$, but independent of $N$, such that
for $k=0,\dots,N-2$ :
\begin{align*}
\summ v(k,m) \leq& C h^{n+1}\lp 1 + {1\over h^{n-p+\undeux}(N-k-1)^{n-p+\undeux}}\rp \\
	\leq& C \lp h^{n+1}+{h^{p+\undeux}\over (N-k-1)^{n-p+\undeux}}\rp.
\end{align*}
Therefore, there exists a constant $C$ as above
such that
\begin{align*}
\summ\sum_{k=0}^{N-2} v(k,m) \leq & C\lp h^n+h^{p+\undeux}\sum_{k=0}^{N-2}
		{1\over(N-k-1)^{n-p+\undeux}}\rp\\
	\leq &C\lp h^n+h^{p+\undeux}\sum_{l=1}^{N-1}
		{1\over l^{n-p+\undeux}}\rp
	\leq C h^{p+\undeux},
\end{align*}
which concludes the proof.
\endproof

\subsection{Decomposition of the weak error.}
We follow the classical decomposition introduced in \cite{TalayTubaro}.
The definition of $u_{\lambda}(t,x)$ in section 3.1 yields
\begin{align*}
E\lv X^N(T)\rvd_\Hp - E\lv X(T)\rvd_\Hp = &
	\summ\lmp\lp E\lv\xnlm(T)\rvd - E\lv\xlm(T)\rvd\rp\\
	=& \summ\lmp\lp E u_\lm\lp T,\xnlm(T)\rp-u_\lm\lp0,\xnlm(0)\rp\rp.
\end{align*}
Let $\dnkm:=\lmp\lp Eu_\lm\lp\tkk,\xnlm\ltkk\rp - Eu_\lm\lp t_k,\xnlm\ltk\rp\rp$;
then
$$E\lv X^N(T)\rvd_\Hp - E\lv X(T)\rvd_\Hp = \summ\sum_{k=0}^{N-1}\dnkm.$$
Note that using Lemmas 3.3, 3.4 and (3.4) we deduce that for any
$k=0,\dots,N-1$
$$E\int_{t_k}^{t_{k+1}}\lv \gamma_{\lambda}^{k,N}(t){\partial u\over\partial x}
(t,X_{\lambda}^N(t))\rvd dt <\infty.$$
From now, we do not justify that the stochastic integral are centered.
It\^o's formula and Lemma \ref{lem-bzX}, we imply that for $k=0,\dots,N-1$

\def\knlm{^{k,N}_\lm}
\begin{align*}
\dnkm=&\lmp E\int_{t_k}^\tkk\lb\dt u_\lm+\bknlm(t)\dx u_\lm
	+\undeux\lv\gknlm(t)\rvd\dxx u_\lm\rb\txnlm dt\\
&=\lmp E\int_\tk^\tkk\lb I\knlm(t)+\undeux J\knlm(t)\rb dt,
\end{align*}
where
\begin{align}
I\knlm(t) :=& \lp\bknlm(t)+\lm\xnlm(t)\rp\dx u_\lm\txnlm,\label{eq-IkN}\\
J\knlm(t) :=& \lp\lv\gknlm(t)\rvd-1\rp\dxx u_\lm\txnlm. \label{eq-JkN}
\end{align}
This yields the following decomposition:
\def\sumk{\sum_{k=0}^{N-1}}
\begin{align}
E\lv X^N(T)\rvd_\Hp - E\lv X(T)\rvd_\Hp =& \summ \delta^N(N-1,m)
	+\summ\sum_{k=0}^{N-2}\lmp E\int_\tk^\tkk I\knlm(t)dt \nonumber\\
	&+\undeux\summ\sum_{k=0}^{N-2}\lmp E\int_\tk^\tkk J\knlm(t)dt.\label{eq-WE}
\end{align}
Now we study each term of this decomposition.

\begin{lemma} \label{lem-deltaNm}
There exists a constant $C$, independant of $N$, such that
$$\summ\lv\delta^N(N-1,m)\rv\leq Ch^{p+\undeux}.$$
\end{lemma}
This study is similar to the third step of \cite{Debussche}, page 97.
\proof
Using the definition of $u_{\lambda_m}(t,x)$
\eqref{eq-u} and \eqref{eq-xnld}, we have
\begin{align*}
u_\lm\lp t_N, X^N_\lm(t_N)\rp =& \lv X^N_\lm(t_N)\rvd
	= {1\over\lp1+\lm h\rp^2}\lv X^N_\lm(t_{N-1})+\Delta W_m(N)\rvd,\\
u_\lm\lp t_{N-1}, X^N_\lm(t_{N-1})\rp =& {1-e^{-2\lm h}\over2\lm} 
	+e^{-2\lm h}\lv X^N_\lm(t_{N-1})\rvd.
\end{align*}
By independence between $\Delta W_m(N)$ and $X^N_\lm(t_{N-1})$, we have
\begin{align*}
\delta^N(N-1,m) =& \lmp\lb {1\over\lp1+\lm h\rp^2}-e^{-2\lm h}\rb
	E\lv\xnlm(t_{N-1})\rvd \\
	&+ {h\over\lambda^p_m\lp1+\lm h\rp^2} -{1-e^{-2\lm h}\over 2\lambda^{1+p}_m}.
\end{align*}
\def\dul{\delta_1\lp\lm\rp}
\def\ddl{\delta_2\lp\lm\rp}
\def\dtl{\delta_3\lp\lm\rp}
Let $\dul := {1-2e^{-2\lm h}\over2\lm^{1+p}}$,
$\ddl := {h\over\lm^p\lp1+\lm h\rp^2}$, and
$$\dtl := \lmp\lb{1\over\lp1+\lm h\rp^2}-e^{-2\lm h}\rb E\lv\xnlm(t_{N-1})\rvd.$$
With these notations we have
$$\delta^N(N-1,m)\leq\dul+\ddl+\dtl.$$
First, we study $\dul$.
Since ${1-e^{-2\lambda h}\over2\lambda} = \int_0^he^{-2\lambda x}dx$, using
Lemma \ref{lem-a1}, we obtain
\begin{align}
\summ\dul =& \int_0^h\summ\lmp e^{-2\lm x}dx 
	\leq C\int_0^h x^{p-\undeux}dx 
	= Ch^{p+\undeux}. \label{eq-dul}
\end{align}
Now we study $\ddl$.
Since $\lp x\in[0,\infty)\mapsto x^{-2p}(1+x^2h)^2\rp$ is decreasing, we have
for $p\in[0,\undeux)$
\begin{align}
\summ\ddl \leq& Ch\int_0^\infty{1\over x^{2p}\lp1+x^2h\rp^2}dx 
	\leq  Ch^{p+\undeux}\int_0^\infty {y^{-2p}\over(1+y^2)^2}dy 
	\leq  Ch^{p+\undeux}. \label{eq-ddl}
\end{align}
Finally, we study $\dtl$.
Using Lemma \ref{lem-EX}, we have
$$\dtl \leq \lmp\lb{1\over(1+\lm h)^2}-e^{-2\lm h}\rb{1\over2\lm}.$$
Since ${1\over\lp1+\lambda h\rp^2}-e^{-2\lambda h}
	=2\lambda\int_0^h\lb e^{-2\lambda x}-{1\over\lp1+\lambda x\rp^3}\rb dx$,
we have
$$\dtl \leq \lmp\int_0^h\lb e^{-2\lm x}+{1\over(1+\lm x)^3}\rb dx.$$
Using Lemma \ref{lem-a1}, we have for $p\in[0,\undeux)$
$$\summ\lmp\int_0^h e^{-2\lm x}dx \leq C\int_0^h x^{p-\undeux}dx
	\leq Ch^{p+\undeux}.$$
Now since for $x\geq0$ the map $\lp y\in\R_+\mapsto y^{-2p}(1+y^2x)^{-3}\rp$ is decreasing, we have
for $p\in[0,\undeux)$
\begin{align*}
\summ{\lmp\over(1+\lm x)^3} \leq C \int_0^\infty {1\over y^{2p}\lp1+y^2x\rp}dy
	\leq C x^{p-\undeux}\int_0^\infty{1\over z^{2p}\lp1+z^2\rp^3}dz 
	\leq Cx^{p-\undeux},
\end{align*}
and hence Fubini's theorem yields
$$\summ\int_0^h{\lmp\over\lp1+\lm x\rp^3}dx \leq C\int_0^h x^{p-\undeux}dx
	\leq Ch^{p+\undeux}.$$
The above inequalities imply $\summ\dtl\leq Ch^{p+\undeux}.$
This inequality, \eqref{eq-dul} and \eqref{eq-ddl} give the stated upper estimate.
\endproof
\begin{lemma}
\label{lem-J}
There exists a constant $C>0$, independent of $N$, such that
$$\summ\sumk\lmp E\int_\tk^\tkk \lv J\knlm(t)\rv dt \leq C h^{p+\undeux}.$$
\end{lemma}
\proof
Using Lemma \ref{lem-bzX}, we have
$$\lv\gknlm(t)\rvd-1 = -{2(t-t_k)\lm\over1+\lm h}
	+{\lv t-t_k\rvd\lm^2\over(1+\lm h)^2}.$$
Using \eqref{eq-uxx} and \eqref{eq-JkN}, we have
$$\lmp E\int_{t_k}^\tkk\lv J^{k,N}_\lm(t)\rv dt 
	\leq C\lp\lambda^{1-p}_mh^2+\lambda^{2-p}_mh^3\rp e^{-2\lm(T-\tkk)}.$$
Lemma \ref{lem-at} concludes the proof.
\endproof
\begin{lemma}\label{lem-I}
There exists a constant $C>0$, independant of $N$, such that
$$\summ\sum_{k=0}^{N-2}\lmp E\int_\tk^\tkk\lv I^{k,N}_\lm(t)\rv dt
	\leq Ch^{p+\undeux}.$$
\end{lemma}
\proof
\def\tkkxnlm{\lp \tkk,\xnlm\ltkk\rp}
\def\Iu{I^{k,N}_{1,\lm}(t)}
\def\Id{I^{k,N}_{2,\lm}(t)}
Let $\Iu := E\bknlm(t)\dx u_\lm\txnlm + E\lm\xnlm\ltkk\dx u_\lm\tkkxnlm$ and
$\Id := -\lm E\xnlm\ltkk\dx u_\lm\tkkxnlm +\lm E\xnlm(t)\dx u_\lm\txnlm$.
Using \eqref{eq-IkN}, we have
\begin{equation}
E I^{k,N}_\lm(t) = \Iu + \Id.
\label{eq-Idec}
\end{equation}

First we study $\Iu$. Using \eqref{eq-ux}, we know that $\dx u_\lm\in C^{1,2}$.
So using It\^o's formula and Lemma \ref{lem-bzX}, we have
\def\sxnlm{\lp s,\xnlm(s)\rp}
\begin{align}
d\dx u_\lm\sxnlm =& \lb \dtx u_\lm+\bknlm(s)\dxx u_\lm\rb\sxnlm ds\nonumber\\
	&+\gknlm(s)\dxx u_\lm\sxnlm dW_\lm(s) \label{eq-duX}
\end{align}
Using this equation, Lemma \ref{lem-bzX} and the It\^o formula we deduce
\begin{align*}
d\lc\bknlm(s)\dx u_\lm\right.&\left.\sxnlm \rc= \lb\bknlm(s)\dtx u_\lm
	+\lv\bknlm(s)\rvd\dxx u_\lm\right. \\ 
	&\left.+z\knlm(s)\gknlm(s)\dxx u_\lm\rb\sxnlm ds\\
	&+\lb\bknlm(s)\gknlm(s)\dxx u_\lm +z\knlm(s)\dx u_\lm \rb
		\sxnlm dW_\lm(s).
\end{align*}
Integrating between $t$ and $\tkk$, taking expectation, 
and using the fact that $\beta^{k,N}_{\lambda_m}(t_{k+1})=-\lambda_m
X^N_{\lambda_m}(t_{k+1})$, so that $I^{k,N}_{1,\lambda_m}(t_{k+1})=0$,
we obtain
\begin{align}
\Iu = -E\int_t^\tkk&\lb\bknlm(s)\dtx u_\lm +\lv\bknlm(s)\rvd\dxx u_\lm\right.\nonumber\\
	&\left.+z\knlm(s)\gknlm(s)\dxx u_\lm\rb\sxnlm ds.
	\label{eq-Iu}
\end{align}
Using \eqref{eq-utx} and Lemma \ref{lem-bx}, we have for $s\in[t,\tkk]$
$$E\bknlm(s)\dtx u_\lm\sxnlm = 4\lm e^{-2\lm(T-s)} E\bknlm(s)\xnlm(s)
	\leq C\lm e^{-2\lm(T-\tkk)},$$
and hence
$$\lmp\int_\tk^\tkk dt\int_t^\tkk ds E\bknlm(s)\dtx u_\lm\sxnlm
	\leq C\lambda^{1-p}_mh^2e^{-\lm(T-\tkk)}.$$
Using Lemma \ref{lem-at}, and the above inequality, we deduce
\begin{equation}
\summ\sum_{k=0}^{N-2}\lmp\int_\tk^\tkk dt\int_t^\tkk ds E\bknlm(s)\dtx u_\lm\sxnlm
	\leq Ch^{p+\undeux}.
	\label{eq-Iu1}
\end{equation}

Using \eqref{eq-uxx} and Lemma \ref{lem-bx}, we have for $s\in[t_k,t_{k+1}]$
$$E\lv\bknlm(s)\rvd\dxx u_\lm\sxnlm = 4\lm e^{-2\lm(T-s)} \leq 4\lm e^{-2\lm(T-\tkk)},$$
so that
$$\lmp\int_\tk^\tkk dt\int_t^\tkk ds E\lv\bknlm(s)\rvd\dxx u_\lm\sxnlm
	\leq C\lambda^{1-p}_mh^2e^{-2\lambda(T-\tkk)}.$$
Thus, Lemma \ref{lem-at} yields
\begin{equation}
\summ\sum_{k=0}^{N-2} \lmp\int_\tk^\tkk dt\int_t^\tkk ds E\lv\bknlm(s)\rvd\dxx u_\lm\sxnlm
	\leq Ch^{p+\undeux}.
	\label{eq-Iu2}
\end{equation}
Using equations \eqref{eq-uxx} and Lemma 3.4 we have for all
$s\in[t,\tkk]$
\begin{align*}
E\lv z\knlm(s)\gknlm(s)\dxx u_\lm\sxnlm\rv =& {2\lm\over1+\lm h}
	\lp1-{(s-\tk)\lm\over1+\lm h}\rp e^{-2\lm(T-s)}\\
	\leq& C\lm e^{-2\lm(T-\tkk)}.
\end{align*}
Therefore, we obtain
$$\lmp\int_\tk^\tkk dt\int_t^\tkk ds E\lv z\knlm(s)\gknlm(s)\dxx u_\lm\sxnlm\rv
	\leq C\lambda^{1-p}_mh^2 e^{-2\lm(T-\tkk)}.$$
Using once more Lemma \ref{lem-at}, we deduce
$$\summ\sum_{k=0}^{N-2}
	\lmp\int_\tk^\tkk dt\int_t^\tkk ds E\lv z\knlm(s)\gknlm(s)\dxx u_\lm\sxnlm\rv
	\leq Ch^{p+\undeux}.$$
Plugging this inequality together with \eqref{eq-Iu1} and \eqref{eq-Iu2} into \eqref{eq-Iu}
gives us
\begin{equation}
\summ\sum_{k=0}^{N-2}\lmp E\int_\tk^\tkk\lv \Iu\rv dt
	\leq Ch^{p+\undeux}.
	\label{eq-Iubound}
\end{equation}

Now we study $\Id$.
Using Lemma \ref{lem-bzX}, equation \eqref{eq-duX} and the It\^o formula
we have
\begin{align*}
d\xnlm(s)\dx& u_\lm\sxnlm = \lb\xnlm(s)\dtx u_\lm
	+\xnlm(s)\bknlm(s)\dxx u_\lm\right.\\
	&\left. +\bknlm(s)\dx u_\lm
	+\lv\gknlm(s)\rvd\dxx u_\lm\rb\sxnlm ds \\
	&+ \lb \gknlm(s)\dx u_\lm +\xnlm(s)\gknlm(s)\dxx u_\lm \rb
		\sxnlm dW_\lm(s)
\end{align*}
So integrating between $t$ and $\tkk$ and taking expectation, we obtain
\begin{align}
\Id = -\lm E\int_t^{\tkk} &\lb \xnlm(s)\dtx u_\lm
	+\bknlm(s)\dx u_\lm +\xnlm(s)\bknlm(s)\dxx u_\lm\right.\nonumber\\
	&\left.+\lv\gknlm(s)\rvd\dxx u_\lm \rb\sxnlm ds.
\label{eq-Id}
\end{align}
Using equation \eqref{eq-utx} and Lemma \ref{lem-bx}, we have for
all $s\in[t,\tkk]$
\def\dsx{{\partial^2\over\partial t\partial x}}
\begin{align*}
\lm E\xnlm(s)\dsx u_\lm\sxnlm
	=& 4\lambda^2_me^{-2\lm(T-s)} E\lv\xnlm(s)\rvd\\
	\leq& C\lambda^2_m ({1\over\lm}+h)e^{-2\lm(T-\tkk)}.
\end{align*}
Therefore,
$$\lmp\int_\tk^\tkk\int_t^\tkk\lm E\xnlm(s)\dsx u_\lm\sxnlm \leq
	C\lp\lambda^{1-p}_mh^2+\lambda^{2-p}_mh^3\rp e^{-2\lm(T-\tkk)},$$
and using Lemma \ref{lem-at}, we deduce
\begin{equation}
\summ\sum_{k=0}^{N-2}
	\lmp\int_\tk^\tkk dt\int_t^\tkk ds\lm E\xnlm(s)\dsx u_\lm\sxnlm \leq
	Ch^{p+\undeux}.
\label{eq-Id1}
\end{equation}
The equation \eqref{eq-ux} and Lemma \ref{lem-bx} yield for all 
$s\in[t,\tkk]$
$$\lm E\bknlm(s)\dx u_\lm\sxnlm
	= 2\lm e^{-2\lm(T-s)}E\bknlm(s)\xnlm(s)
	\leq C\lm e^{-2\lm(T-\tkk)}.$$
This upper estimate implies
$$\lmp\int_\tk^\tkk dt\int_t^\tkk ds\lm E\bknlm(s)\dx u_\lm\sxnlm
	\leq C\lambda^{1-p}_mh^2e^{-2\lm(T-\tkk)},$$
and Lemma \ref{lem-at} yields
\begin{equation}
\summ\sum_{k=0}^{N-2}
\lmp\int_\tk^\tkk dt\int_t^\tkk ds\lm E\bknlm(s)\dx u_\lm\sxnlm
	\leq Ch^{p+\undeux}.
\label{eq-Id2}
\end{equation}
Using equation \eqref{eq-uxx} and Lemma \ref{lem-bx}, we have for all
$s\in[t,\tkk]$
$$\lm E\xnlm(s)\bknlm(s)\dxx u_\lm\sxnlm
	\leq C\lm e^{-2\lm(T-\tkk)}.$$
Therefore, we obtain
$$\lmp\int_{t_k}^\tkk dt\int_t^\tkk ds\lm E\xnlm(s)\bknlm(s)\dxx u_\lm\sxnlm
	\leq C\lambda^{1-p}_mh^2e^{-2\lm(T-\tkk)},$$
and Lemma \ref{lem-at} implies
\begin{equation}
\summ\sum_{k=0}^{N-2}
\lmp\int_{t_k}^\tkk dt\int_t^\tkk ds\lm E\xnlm(s)\bknlm(s)\dxx u_\lm\sxnlm
	\leq Ch^{p+\undeux}.
\label{eq-Id3}
\end{equation}
Finally, \eqref{eq-uxx} and Lemma \ref{lem-bzX} imply that for all $s\in[t,\tkk]$
$$\lm E\lv\gknlm(s)\rvd\dxx u_\lm\sxnlm
	\leq C\lm e^{-2\lm(T-\tkk)}.$$
This yields
$$\lmp\int_\tk^\tkk dt\int_t^\tkk ds\lm E\lv\gknlm(s)\rvd\dxx u_\lm\sxnlm
	\leq C\lambda^{1-p}_mh^2e^{-2\lm(T-\tkk)},$$
and Lemma \ref{lem-at} implies
$$\summ\sum_{k=0}^{N-2}
	\lmp\int_\tk^\tkk dt\int_t^\tkk ds\lm E\lv\gknlm(s)\rvd\dxx u_\lm\sxnlm
	\leq Ch^{p+\undeux}.$$
Plugging this inequality together with \eqref{eq-Id1} -
\eqref{eq-Id3} into \eqref{eq-Id}, we deduce
$$\summ\sum_{k=0}^{N-2}\lmp E\int_\tk^\tkk\lv\Id\rv dt
	\leq Ch^{p+\undeux}.$$
This equation together with \eqref{eq-Idec} and \eqref{eq-Iubound} conclude
the proof.
\endproof
Theorem 2.1 is a straightforward consequence of equation \eqref{eq-WE}
and Lemmas \ref{lem-deltaNm}-\ref{lem-I}.
\subsection*{Acknowledgments:}
The author wishes to thank Annie Millet for many helpful comments.


\begin{thebibliography}{99}
\bibitem{Aboura}
Aboura O.,
Weak error expansion of the implicit Euler scheme
\bibitem{AL}
Andersson A., Larsson S.,
Weak convergence for a spatial approximation of the nonlinear stochastic heat equation,
arixv 1212.5564
\bibitem{BallyTalay}
Bally V., Talay D.,
The law of the Euler scheme for stochastic differential equations:
I. Convergence rate of the distribution function. 
\textit{Prob. Th. Rel. Fields}, 104-1:43--60, 1996.
\bibitem{DaPrato}
Da Prato G., Zabczyk J.,
{\it Stochastic equations in infinite dimensions}
\bibitem{deBouard}
de Bouard A., Debussche A.,
Weak and strong order of convergence of a semi discrete
scheme for the stochastic Nonlinear Schrodinger equation, 
\textit{Applied Mathematics and
    Optimization journal}, 54, 3, pp. 369--399 2006.
\bibitem{Debussche}
Debussche A.,
Weak approximation of stochastic partial differential equations: the
nonlinear case, \textit{Math of Comp}, 80, pp. 89--117 2011.
\bibitem{Printems}
Debussche A., Printems J.,
Weak order for the discretization of the stochastic heat
equation, \textit{Math of Comp}, 78, 266, pp. 845--863 2009.
\bibitem{G1}
Gy\"ongy I.,
Lattice approximations for stochastic quasi-linear parabolic partial
differential equations driven by space-time white noise. I
{\it Potential Anal.} 9, no. 1, 1--25, 1998.
\bibitem{G2}
Gy\"ongy I.,
Lattice approximations for stochastic quasi-linear parabolic partial
differential equations driven by space-time white noise. II
{\it Potential Anal.} 11, no. 1, 1--37, 1999.
\bibitem{Hausenblas}
Hausenblas E., Weak approximation for semilinear stochastic
evolution equations,
{\it Stochastic Analysis and Related Topics VIII}, 2003
\bibitem{Kloeden}
Kloeden P. E., Platen E.,
{\it Numerical solution of stochastic differential equations.}
Applications of Mathematics (New York), 23.
{\it Springer-Verlag, Berlin}, 1992.
\bibitem{KLL}
Kov\'acs M., Larsson S., Lindgren F.,
Weak convergence of finite element approximations of linear 
stochastic evolution equations with additive noise,
{\it BIT} {\bf 52}, 2012.
\bibitem{KLL2}
Kov\'acs M., Larsson S., Lindgren F.,
Weak convergence of finite element approximations of linear stochastic evolution equations with additive noise II. Fully discrete schemes,
{\it BIT Numer. Math.}, 2012.  
\bibitem{Nualart}
Nualart D.,
{\it The Malliavin Calculus and Related Topics. Second Edition.}
Probability and its Applications (New York).
{\it Springer-Verlag, Berlin,} 2006.
\bibitem{PJ}
Printems J.,
On the discretization in time of parabolic stochastic partial differential
equations,
{\it Math. Model. and Numer. Anal.}, 35(6), 1055-1078, 2001.

\bibitem{TalayTubaro}
Talay D., Tubaro L.,
Expansion of the global error for numerical schemes solving SDEs
{\it Stoch. Anal. and App.}, 8(4),483--509, 1990.
\end{thebibliography}
\end{document}